\documentclass[12pt]{amsart} \usepackage{amssymb}
\usepackage{amsfonts} \usepackage{latexsym}
\usepackage{amscd}
\newcounter{cs} \stepcounter{cs} \newcounter{ds} \stepcounter{ds}

\newcommand{\casos}{\begin{itemize}}
\newcommand{\fcasos}{\end{itemize}\setcounter{cs}{1}}

\newcommand{\mul}{{\mathcal M}(R)}
\newcommand{\ra}{\rightarrow}

\newcommand{\ilim}{\varprojlim}
\newcommand{\ol}{\overline}

\newcommand{\nfi}{\varphi}
\newcommand{\bsr}{\operatorname{bsr}}
\newcommand{\id}{\operatorname{id}}

\newtheorem{lem}{Lemma}[section]
\newtheorem{corol}[lem]{Corollary}
\newtheorem{theor}[lem]{Theorem}
\newtheorem{prop}[lem]{Proposition}

\newtheorem{defi}[lem]{Definition}
\newtheorem{exem}[lem]{Examples}

\begin{document}

\title{Inverse Limits of Rings and Multiplier Rings}

\author{Gert K. Pedersen \& Francesc Perera}

\address{Department of Mathematics, University of Copenhagen,
  Universitetsparken 5, DK-2100, Copenhagen \O}
\address{Department of Pure Mathematics, Queens University Belfast,
BT7 1NN, Northern Ireland}

\email{gkped@math.ku.dk \& perera@qub.ac.uk}

\date{July 2003}

\thanks{Research supported by the Danish Research Council, SNF, and the
Nuffield Foundation}

\keywords{Bass rank, Exchange ring, Multiplier ring, $QB-$ring}


\begin{abstract} It is proved that the exchange property, the
Bass stable rank and the quasi-Bass property are all preserved
under surjective inverse limits. This is then applied to
multiplier rings by showing that these in many cases can be
obtained as inverse limits.
\end{abstract}

\maketitle

\section*{Introduction}

Given a sequence $(R_n)$ of rings with connecting morphisms (i.e.
ring homomorphisms) $\pi_n\colon R_n \ra R_{n-1}$ for all $n$
(setting $R_0=R_1$ and $\pi_1=\mathrm{id}$) we define the {\em
inverse limit} as the ring $\ilim R_n$ of strings $x=(x_n)$ in
$\prod R_n$, i.e. sequences such that $\pi_n(x_n)=x_{n-1}$ for all
$n$. For each $m$ there is a natural morphism $\rho_m\colon \ilim
R_n \ra R_m$ (the {\em coordinate evaluation}) obtained by
evaluating a string $x=(x_n)$ at $m$, and we see that
$\pi_n\circ\rho_n=\rho_{n-1}$ for every $n$. The ring $R=\ilim
R_n$ has the universal property that for each coherent sequence of
morphisms $\sigma_n\colon S\ra R_n$ from a ring $S$ (i.e.
$\pi_n\circ\sigma_n =\sigma_{n-1}$ for all $n$) there is a unique
morphism $\sigma\colon S\ra R$ such that $\sigma_n
=\rho_n\circ\sigma$ for all $n$.

If for each $m$ we let $S_m=\rho_m(\ilim R_n) \subset R_m$, then
$\pi_n(S_n)=S_{n-1}$ and $\ilim R_n = \ilim S_n$. This shows that
every ring which can be obtained as an inverse limit can also be
obtained as an inverse limit in which each morphism $\pi_n$ is
surjective. We shall refer to this case as a {\em surjective
inverse limit}, and will concentrate exclusively on it.

By their very construction inverse limits tend to be large. For
example, if the morphisms are not eventually constant, an inverse
limit will be uncountable even if the rings are finite. This might
be considered a detracting factor. In this paper we shall try to
redeem the construction by showing that in many instances
multiplier rings (which we expect to be large) can be obtained as
inverse limits. Since the structural properties of inverse limits
are good, we hereby obtain information about multiplier rings that
would otherwise seem unreachable.

We prove in Section 1 that a surjective inverse limit of exchange
rings is again an exchange ring. In Section 2 we show that the
Bass stable rank of a surjective inverse limit of rings $(R_n)$ is
the supremum of the Bass stable ranks of the $R_n$'s. In Section 3
we extend this to an important infinite case by proving that the
surjective inverse limit of quasi-Bass rings is again a quasi-Bass
ring.

In Sections 4 and 5 we establish the basic properties of
approximate units and multiplier rings that we shall need. In
particular we show that for every proper morphism $\nfi\colon
R\rightarrow S$ between non-degenerate rings there is a unique
extension $\ol\nfi\colon \mul\rightarrow \mathcal M (S)$ between
their multiplier rings. We also show that $\ol\nfi$ is strictly
continuous, cf.~\cite{arap}. If $R$ is $\sigma-$unital, i.e. has a
countable approximate unit, and $\nfi$ is surjective we show that
$\ol\nfi$ is surjective as well, thereby giving an algebraic
analogue of the Tietze extension theorem (see e.g.~\cite{eng}).

Finally in Section 6 we combine our results by showing that many
multiplier rings can be obtained as inverse limits. Thus if
$(R_n)$ is a sequence of $\sigma-$unital rings with surjective
morphisms $\pi_n\colon R_n \rightarrow R_{n-1}$ we prove that
$\mathcal M(\ilim R_n) = \varprojlim\mathcal M (R_n)$. Going
further we show that if a semi-prime $\sigma-$unital ring $R$ has
two sequences of ideals $(I_n)$ and $(J_n)$, one increasing to
$R$, the other decreasing to $0$, such that $I_n\cap J_n =0$ and
$R/J_n$ is unital for every $n$, then $\mul=\ilim R/J_n$.

The inspiration for this paper comes from~\cite{bp5}, where a
similar list of results are obtained for $C^*-$algebras. In this
category one must of course consider only {\em bounded strings} as
elements in the inverse limit, so although our results are the
same as in~\cite{bp5}, our proofs are completely different -
sometimes harder, sometimes easier. It is worth mentioning that
inverse limits of $C^*-$algebras in the algebraic sense (with
arbitrary, unbounded strings) have been explored by Phillips
in~\cite{phillips} and~\cite{phillips2} to obtain examples of {\em
pro $C^*-$algebras}. The aim was to create a non-commutative
analogue of normal spaces that are not necessarily locally
compact. The unpublished notes ~\cite{good2} by Goodearl have also
influenced our work, although we fail to answer the question that
motivated them: Will a surjective inverse limit of separative
exchange rings itself be separative? A negative answer would
provide a solution to the Fundamental Separativity Problem
(see~\cite{4au}).

\vskip1cm

\section{Exchange Rings}

\bigskip

Recall from~\cite{goodwar} and~\cite{nic} that a (unital) ring $R$
is an {\em exchange ring} if for each $x$ in $R$ there is an
idempotent $e$ in $xR$ such that $1-x=(1-e)(1-y)$ for some $y$ in
$R$. This is not the original definition, which concerns a finite
exchange property for $R-$modules, see~\cite{war}, but an
equivalent formulation better suited for our purposes.

The class of exchange rings is pleasantly large and includes all
(von Neumann) regular rings, all $\pi-$regular rings, the
semi-perfect rings (identified with the semi-local exchange rings)
and the $C^*-$algebras of real rank zero (identified with the
exchange $C^*-$algebras by~\cite[Theorem 7.2]{4au}).

We shall show that every surjective inverse limit of exchange
rings again has the exchange property. As an hors d'\oe{}uvre we
present the following simpler result from~\cite{good2} (see
also~\cite[Example 1.10]{good}). Recall that an element $x$ in a
ring $R$ is said to be \emph{regular} provided that $x=xyx$ for
some $y$ in $R$. We shall refer to such an element $y$ as a
\emph{partial inverse} for $x$. If all elements of $R$ are regular
we say that $R$ is \emph{(von Neumann) regular}.

\begin{prop} {\rm (Goodearl)} If $R$ is the surjective
inverse limit of a sequence $(R_n)$ of regular rings, then $R$ is
also regular.
\end{prop}

\begin{proof} Put $R= \varprojlim R_n$ and let
$\pi_n\colon R_n\rightarrow R_{n-1}$ denote the connecting
morphisms. If $x=(x_n)$ is an element in $R$ and if for all $k<n$
we have found elements $y_k$ in $R_k$ such that
$\pi_k(y_k)=y_{k-1}$ and $x_ky_kx_k=x_k$, we choose an element $y$
in $R_n$ with $\pi_n(y)=y_{n-1}$. Then $u=x_n-x_nyx_n \in
\ker\pi_n$, and since this ideal is a regular ring in its own
right we can find a $v$ in $\ker\pi_n$ such that $uvu=u$. By
computation
\begin{align*}
&x_n-x_nyx_n=(x_n-x_nyx_n)v(x_n-x_nyx_n)\\
=\;&x_nvx_n-x_nvx_nyx_n-x_nyx_nvx_n+x_nyx_nvx_nyx_n\,,
\end{align*}
and it follows that the element $y_n=y+v-yx_nv-vx_ny+yx_nvx_ny_n$
is a partial inverse for $x_n$ with $\pi_n(y_n)=y_{n-1}$.

By induction we can therefore find an element $y=(y_n)$ in $R$
such that $x=xyx$, so every element in $R$ is regular, as desired.
\end{proof}

\bigskip

Note that in the definition of an exchange ring the role of the
unit is superfluous. Thus we say that a (not necessarily unital)
ring $R$ is an exchange ring if for each $x$ in $R$ there is an
idempotent $e$ in $xR$ such that $x=e+y-ey$ for some $y$ in $R$.
This idea was successfully exploited by Ara in~\cite{pext}. If $I$
is a non-unital exchange ring embedded as a two-sided ideal of a
unital ring $R$, we will adopt the terminology from~\cite{pext}
and say that $I$ is an \emph{exchange ideal of} $R$.

\bigskip

\begin{lem}\label{aral}
Let $I$ be a two-sided exchange ideal in a unital ring $R$. If $x$
and $y$ are elements in $R$ and $p$ is an idempotent such that
$p-xy\in I$, there is an idempotent $q$ in $I$ and an element $r$
in $pRp$ with $p-r$ in $I$, such that both elements $a=(p-q)x$ and
$b=yr$ are regular and partial inverses for one another. By
construction $a-px\in I$ and $b-yp\in I$, so $ab-p\in I$.
\end{lem}

\bigskip

This is a restatement of~\cite[Lemma 2.1]{pext}. The result is
obtained by applying the exchange condition to the element
$p-pxyp$ in $pIp$ to get $q$.

If $R$ is a non-unital ring, it is sometimes convenient to adjoin
a unit in order to obtain the ring $R^+=R\oplus\mathbb Z$, with
elementwise addition and multiplication given by the rule
$(x,n)(y,m)=(xy+mx+ny,nm)$. In this way $R$ sits as a two-sided
ideal in $R^+$ with quotient $\mathbb Z$.

\bigskip

\begin{lem}\label{lift3} Let $\pi\colon R \to S$ be a
surjective morphism between (not necessarily unital) exchange
rings, and let $\ol{x}, \ol{y}$ and $\ol{z}$ be elements in $S$
such that with $\ol{e}=\ol{x}\,\ol{y}$ we have
\[
\ol{y}=\ol{y}\,\ol{e},\quad\quad 1-\ol{e}=(1-\ol{x})(1-\ol{z}),
\quad\quad (1-\ol{z})\ol{e}=0\,.
\]
For each choice of $x$ in $R$ with $\pi(x)=\ol{x}$ there are then
elements $y$ and  $z$ in $R$  with $\pi(y)=\ol{y}$ and
$\pi(z)=\ol{z}$, such that with $e=xy$ we have
\[
y=ye,\quad\quad 1-e=(1-x)(1-z),\quad\quad (1-z)e=0\,.
\]
\end{lem}

\begin{proof}
This proof amalgamates arguments from~\cite[Theorem 2.2]{pext}
and~\cite[Proposition 1.1]{nic}, which are included for
convenience. As usual we shall write $p\le q$ for idempotents $p$
and $q$ such that $pq=qp=p$. Moreover, if $R$ is not unital we
adjoin a unit to obtain the unital ring $R^+ =R\oplus \mathbb Z$
and we identify $S$ with an exchange ideal in $S^+$, setting
$\pi(1)=1$.

Observe first that $\ol{e}$ is an idempotent by necessity, since
$\ol{e}\,\ol{e}=\ol{x}\,\ol{y}\,\ol{e}=\ol{x}\,\ol{y}=\ol{e}$.
Thus also $\ol{p}_1=\ol{y}\,\ol{x}$ is an idempotent. Since $R$ is
an exchange ring $\ol{p}_1$ can be lifted to an idempotent $p_1$
in $R$ by~\cite[Theorem 2.2]{pext}. Let $y_1$ and $z_1$ be any
lifts in $R$ of the elements $\ol{y}$ and $\ol{z}$. By
Lemma~\ref{aral} (with $I=\ker(\pi)$ and exchanging the roles of
$x$ and $y$) we can find regular elements $a_1$ and $b_1$ in $R$
that are partial inverses for one another, such that
$a_1=(p_1-q_1)y_1$ and $b_1=xr_1$ for some idempotent $q_1$ in $I$
and some $r_1$ in $p_1Rp_1$ with $\pi(r_1)=\ol{p}_1$. Thus
$e_1=b_1a_1$ is an idempotent with
\[
\pi(e_1)=\ol{x}\,\ol{p}_1\ol{p}_1\ol{y}=\ol{e}\,\ol{e}\,\ol{e}=\ol{e}\,.
\]

Put $z_2=1-(1-z_1)(1-e_1)$, so that
$\pi(z_2)=1-(1-\ol{z})(1-\ol{e}) =\ol{z}$. Observe that
$\ol{p}_2=(1-\ol{z})(1-\ol{x})$ is an idempotent in $1-S$ and lift
it to an idempotent $p_2$ in $1-R$. Again by Lemma~\ref{aral}
there are regular elements $a_2$ and $b_2$ in $1-R$, partial
inverses for one another, such that $a_2=(p_2-q_2)(1-z_2)$ and
$b_2=(1-x)r_2$ for some idempotent $q_2$ in $I$ and some element
$r_2$ in $p_2R^+ p_2$ with
$\pi(r_2)=\ol{p}_2=(1-\ol{z})(1-\ol{x})$ (so that $r_2 \in 1-R$ as
well).

The idempotent $f=b_2a_2$ belongs to $(1-x)(1-R)(1-z_2)$ and
satisfies
\[
\pi(f)=(1-\ol{x})\ol{p}_2(1-\ol{z})
=(1-\ol{x})(1-\ol{z})(1-\ol{x})(1-\ol{z})=1-\ol{e}\,.
\]
Since $f=f(1-e_1)$ by construction of $z_2$, we may consider the
idempotent $(1-e_1)f$, where now $(1-e_1)f \le 1-e_1$. Thus
$q=(1-e_1)(1-f)$ is an idempotent in $I$ with $q\le 1-e_1$.

Now apply the exchange property to the element $qxq$ in the corner
ring $qRq = qR^+q \, (\subset I)$ to find an idempotent $t$ in
$qRq$ such that $t=qxqcq$ and $q-t=q(1-x)qdq$ for some elements
$c, d$ in $R$. Put $s=xqct$, and note that $sq=s$ and $qs=t$,
whence $s^2=s$. Since $t\le q\le 1-e_1$, the element
$e_2=e_1+(1-e_1)s$ is an idempotent in $R$ with $e_1 \le e_2$ and
$\pi(e_2)=\ol{e}$. Let $y_2=r_1a_1(1-s)+qct$. Then
$\pi(y_2)=\ol{p}_1\ol{y}=\ol{y}\,\ol{x}\,\ol{y}=\ol{y}$. Moreover,
\begin{equation}\label{xy_2}
xy_2=b_1a_1(1-s)+s=e_1(1-s)+s=e_1+(1-e_1)s=e_2\,.
\end{equation}

The argument in~\cite[Theorem 2.2]{pext} proceeds to show that the
right ideal $A = e_2 R^+ +(1-x)R^+$ equals $R^+$, and does so by
showing that both $1-q$ and $q$ belong to $A$. All we need is the
last assertion, but there seems no easy way to obtain the specific
decomposition of $q$ except by the full argument.

Evidently $e_2e_1\in A$ and $e_2(1-e_1)\in A$. We know that
$e_2e_1=e_1$, and therefore $e_2(1-e_1)=e_2-e_1=(1-e_1)s= s-e_1s$.
Consequently both $e_1\in A$ and $s\in A$. Now $f=(1-x)r_2a_2\in
A$, so also $(1-e_1)f=f-e_1f\in A$. It follows that
$1-q=e_1+(1-e_1)f\in A$. As $s\in A$ and $t=qs$ we now conclude
that $t=s-(1-q)s\in A$. Since
\[
q-t=q(1-x)qdq=(1-x)qdq-(1-q)(1-x)qdq\in (1-x)R+A=A\,,
\]
we finally see that $q=q-t+t\in A$. We can therefore write
\begin{equation}\label{q}
q=e_2u+(1-x)v
\end{equation}
for some elements $u$ and $v$ in $Rq$. In particular, both
elements belong to $I$.

Put $z_3=1-(r_2a_2 +v)$ and $e_3=e_1(1-f)+u$, and note that
$\pi(1-z_3)=\ol{p}_2\ol{p}_2(1-\ol{z})
=(1-\ol{z})(1-\ol{x})(1-\ol{z})=1-\ol{z}$ and
$\pi(e_3)=\ol{e}\,\ol{e}=\ol{e}$. Moreover, by (\ref{q})
\begin{equation}\label{trick}
\begin{split}
(1-x)(1-z_3)+e_2e_3=\;&b_2a_2+(1-x)v+e_2e_1(1-f)+e_2u\\
=\;&f+e_1(1-f)+q=f+e_1(1-f)+(1-e_1)(1-f)=1\,.
\end{split}
\end{equation}

We can now make our final choices as follows:
\[
e=e_2+e_2e_3(1-e_2),\quad\quad y=y_2e,\quad\quad
z=1-(1-z_3)(1-e_2)(1-e)\,.
\]
It is easy to check that
\[
\pi(e)=\ol{e},\quad\quad \pi(y)=\ol{y}\,\ol{e}=\ol{y},\quad\quad
\pi(z) =1-(1-\ol{z})(1-\ol{e})(1-\ol{e})=\ol{z}.
\]
Moreover, by (\ref{xy_2}) and (\ref{trick}) we have the desired
relations:
\begin{align*}
xy &=xy_2e=e_2e=e\\
ye &=y\quad\text{and}\quad (1-z)e=0\\
(1-x)(1-z) &=(1-x)(1-z_3)(1-e_2)(1-e)
=(1-e_2e_3)(1-e_2)(1-e)\\
&=(1-e_2-e_2e_3(1-e_2))(1-e)=(1-e)(1-e)=1-e\,.
\end{align*}
\end{proof}

\bigskip

\begin{theor}\label{exchange}
If  $R$ is the surjective inverse limit of a sequence $(R_n)$ of
exchange rings, then $R$ is also an exchange ring.
\end{theor}

\begin{proof}
Consider an element $x$ in $R$, identified with a string $(x_n)$
in $\prod R_n$. We must then find an idempotent $e=(e_n)$ in $R$
such that $e\in xR$ and $1-e \in (1-x)(1-R)$.

Since $R_1$ is an exchange ring we can find an idempotent $e_1$
and elements $y_1, z_1$ such that $e_1=x_1y_1$ and
$1-e_1=(1-x_1)(1-z_1)$. Evidently we may also assume that
$y_1=y_1e_1$ and $(1-z_1)e_1=0$.

Assume now that for some $n$ we have found elements $y_k$ and
$z_k$ in $R_k$ for $0 \le k\le n$, such that $\pi_k(y_k)=y_{k-1}$
and $\pi_k(z_k)=z_{k-1}$ for all $k\ge 2$, and moreover with
$e_k=x_ky_k$ we have the relations
\[
y_k=y_ke_k,\quad\quad 1-e_k=(1-x_k)(1-z_k),\quad\quad
(1-z_k)e_k=0,
\]
for all $k$. By Lemma~\ref{lift3} we can then find elements
$y_{n+1}$ and $z_{n+1}$ in $R_{n+1}$ with $\pi_{n+1}(y_{n+1})=y_n$
and $\pi_{n+1}(z_{n+1})=z_n$, such that with
$e_{n+1}=x_{n+1}y_{n+1}$ we have
\[
y_{n+1}=y_{n+1}e_{n+1},\quad\quad
1-e_{n+1}=(1-x_{n+1})(1-z_{n+1}), \quad\quad (1-z_{n+1})e_{n+1}=0.
\]
By induction this defines elements $y=(y_n)$ and $z=(z_n)$ in $R$,
such that the element $e=xy$ is an idempotent and
$1-e=(1-x)(1-z)$, as desired.
\end{proof}

\vskip1cm

\section{Bass Stable Rank}

\bigskip

Let $R$ be a unital ring. As usual we say that a row ${\bf
a}=(a_1, \ldots, a_d)$ in $R^{d}$ is \emph{right unimodular}
provided that $a_1R+\cdots+a_dR=R$. To facilitate the computations
with rows in $R^d$ we introduce the $R-$valued inner product ${\bf
a}\cdot{\bf b} =\sum a_ib_i$, so that ${\bf a}$ is right
unimodular precisely if ${\bf a}\cdot{\bf b}=1$ for some ${\bf b}$
in $R^d$. Also, if ${\bf a}=(a_1,\ldots,a_d)\in R^d$ and $s\in R$,
we write $s{\bf a}=(sa_1,\ldots,sa_d)$.

Recall from~\cite{bass} (see also
e.g.~\cite{bassref},~\cite{lam2}) that the {\em Bass stable rank}
of a unital ring $R$ is the smallest number $\bsr (R)$ such that
for each $d\ge\bsr(R)$ every right unimodular row ${\bf a}=(a_0,
\ldots, a_d)$ in $R^{d+1}$ can be reduced to a right unimodular
row in $R^d$ of the form ${\bf a}^r+a_0{\bf b}$ for a suitable row
${\bf b}$ in $R^d$. Here ${\bf a}^r= (a_1, \ldots, a_n)$ and we
regard $R^d$ as a two-sided module over $R$. Evidently this
definition favours right unimodular rows and should be called the
right Bass rank, but it turns out that the analogous concept for
left unimodular rows gives the same lower bound
(~\cite{vas2},~\cite{bassref}). In particular $\bsr(R)=1$ (in
which case we say that $R$ is a {\em Bass ring}) if every equation
$ax+b=1$ implies that $a+by$ is invertible in $R$ for a suitable
$y$, cf.~\cite{vas}.

\bigskip

\begin{lem}\label{jarl}
Let $R$ be a unital ring with $\mathrm{bsr}(R)\leq d$. Given rows
${\bf a}$ and ${\bf x}$ in $R^d$ with ${\bf a}\cdot{\bf x}=1-s$
for some $s$ in $R$, there exist rows ${\bf b}$ and ${\bf y}$ in
$R^d$ and $z$ in $R$ such that
\[
({\bf a}+s{\bf b})\cdot({\bf x}+{\bf y}sz)=1\,.
\]
\end{lem}

\begin{proof}
Since $\mathrm{bsr}(R)\leq d$, the equation ${\bf a}\cdot{\bf x}
+s1=1$ produces rows ${\bf b}$ and ${\bf y}$ in $R^d$ such that
$({\bf a}+ s{\bf b})\cdot {\bf y} = 1$. Define the elements $y
={\bf a}\cdot{\bf y},\, b={\bf b}\cdot{\bf y}$ and $z=1-{\bf
b}\cdot{\bf x}$, and note that $y+sb=1$. By computation we
therefore get
\begin{align*}
&({\bf a}+s{\bf b})\cdot({\bf x}+{\bf y}sz) =1-s+ysz+s(1-z)+sbsz\\
=\;& 1-s+(y+sb)sz+s(1-z)=1-s+sz+s(1-z)=1\,,
\end{align*}
as desired.
\end{proof}

\bigskip

\begin{lem}\label{lift2}
Consider a surjective morphism $\pi\colon R\rightarrow S$ between
unital rings, where $\mathrm{bsr}(R)\leq d$. Assume that ${\bf a}$
and ${\bf x}$ are unital rows in $R^{d+1}$ such that ${\bf
a}\cdot{\bf x}=1$, and that we have chosen rows $\ol{{\bf b}}$ and
$\ol{{\bf y}}$ in $S^d$ such that $(\pi({\bf
a}^r)+\pi(a_0)\ol{{\bf b}})\cdot \ol{{\bf y}}=1$\,. We can then
find ${\bf b}$ and ${\bf y}$ in $R^d$ such that $\pi({\bf
b})=\ol{{\bf b}}$ and $\pi({\bf y})=\ol{{\bf y}}$, and moreover
$({\bf a}^r+a_0{\bf b})\cdot {\bf y}=1$\,.
\end{lem}

\begin{proof}
Take any lifts ${\bf b}'$ of $\ol{{\bf b}}$ and ${\bf y}'$ of
$\ol{{\bf y}}$. Then $({\bf a}^r+a_0{\bf b}')\cdot {\bf y}'=1-t$
for some $t$ in $I=\ker\pi$.  Since $\mathrm{bsr}(R)\leq d$ we can
use Lemma~\ref{jarl} to find ${\bf s}$ and ${\bf t}$ in $I^d$ such
that
\begin{equation}\label{longsum1}
({\bf a}^r +a_0{\bf b}'+{\bf s})\cdot({\bf y}'+{\bf t})=1\,.
\end{equation}
>From the original condition, setting
$b=x_0 - {\bf b}'\cdot{\bf x}^r$, we get
\[
1 ={\bf a}^r\cdot{\bf x}^r+ a_0x_0 = ({\bf a}^r+a_0{\bf b}')\cdot
{\bf x}^r +a_0b\,.
\]
It follows from (\ref{longsum1}) that with $s={\bf s}\cdot({\bf
y}'+{\bf t})$ in $I$ we have
\begin{equation}\label{longsum2}
\begin{split}
1 &= ({\bf a}^r+a_0{\bf b}')\cdot({\bf y}'+{\bf t}) + s\\
&=({\bf a}^r+a_0{\bf b}')\cdot({\bf y}'+{\bf t})
+ ({\bf a}^r+a_0{\bf b}')\cdot {\bf x}^rs +a_0bs\\
&= ({\bf a}^r+a_0{\bf b}')\cdot({\bf y}'+{\bf t}+{\bf x}^rs)
+a_0bs.
\end{split}
\end{equation}
Using Lemma~\ref{jarl} on (\ref{longsum2}) we find rows ${\bf u}$
in $R^d$ and ${\bf v}$ in $(Ra_0bsR)^d \subset I^d$ such that
\[
1=({\bf a}^r+a_0{\bf b}'+a_0bs{\bf u})\cdot({\bf y}'+{\bf t}+{\bf
v})\,.
\]
Clearly the rows ${\bf b}={\bf b}'+bs{\bf u}$ and ${\bf y}={\bf
y}'+{\bf t} +{\bf v}$ verify the desired conditions.
\end{proof}

\bigskip

\begin{theor}\label{bsr}
If $R$ is the surjective inverse limit of a sequence $(R_n)$ of
rings with $\bsr(R_n)\leq d$ for all $n$, then also $\bsr (R) \leq
d$.
\end{theor}

\begin{proof} Let $\pi_n\colon R_n \rightarrow R_{n-1}$ with
$R_0=R_1$ and $\pi_1=\mathrm{id}$, and consider an equation ${\bf
a}\cdot{\bf x}=1$ in $R^{d+1}$. Identify the rows ${\bf a}$ and
${\bf x}$ with strings $({\bf a}_n)$ and $({\bf x}_n)$ in $\prod
R_n^{d+1}$. Write ${\bf a}_n=(a_{n,0},\ldots,a_{n,d})$. Since
$\mathrm{bsr}(R_1)\le d$ there are rows ${\bf b}_1$ and ${\bf
y}_1$ in $R_1^d$ such that $({\bf a}_1^r + a_{1,0}{\bf b}_1)\cdot
{\bf y}_1=1$. Assume now for some $n$ that we have found rows
${\bf b}_k$ and ${\bf y}_k$ in $R_k^d$ for $1\le k\le n$ such that
$({\bf a}_k^r + a_{k,0}{\bf b}_k)\cdot{\bf y}_k=1$, and moreover
$\pi_k({\bf b}_k)={\bf b}_{k-1}$ and $\pi_k({\bf y}_k)={\bf
y}_{k-1}$ for all $k$. By Lemma~\ref{lift2} there are rows ${\bf
b}_{n+1}$ and ${\bf y}_{n+1}$ in $R_{n+1}^d$ such that $({\bf
a}_{n+1}^r + a_{n+1,0}{\bf b}_{n+1})\cdot{\bf y}_{n+1}=1$, and
moreover $\pi_{n+1}({\bf b}_{n+1})={\bf b}_n$ and $\pi_{n+1}({\bf
y}_{n+1})={\bf y}_n$. By induction we can then define rows ${\bf
b}=({\bf b}_n)$ and ${\bf y}=({\bf y}_n)$ in $R^d$ such that
$({\bf a}^r + a_0{\bf b})\cdot{\bf y}=1$ in $R$, which proves that
$\bsr(R)\le d$, as desired.
\end{proof}

\bigskip

\begin{corol}
Let $R=\varprojlim R_n$, where all the morphisms $R_n\rightarrow
R_{n-1}$ are surjective. Then
\[
\mathrm{bsr}(R)=\sup_n\mathrm{bsr}(R_n)\,.
\]
\end{corol}

\begin{proof} For each $n$ the coordinate projections
  $\rho_n \colon R \rightarrow R_n$ are surjective,
whence $\mathrm{bsr}(R_n)\leq \mathrm{bsr}(R)$ (see,
e.g.~\cite[Theorem 4]{vas2}). Thus $\sup_n\mathrm{bsr}(R_n)\leq
\mathrm{bsr}(R)$.

For the converse inequality, we may of course assume that
$\mathrm{bsr}(R_n)\leq d$ for all $n$ (and some finite $d$). By
Theorem~\ref{bsr} $\mathrm{bsr}(R)\leq d$ as well.
\end{proof}

\vskip1cm

\section{$QB-$Rings}

\bigskip

As shown in~\cite{app1} a $QB-$ring is an infinite version of a
Bass ring, i.e. a ring $R$ with $\bsr(R)=1$. To arrive at the
definition we replace the set $R^{-1}$ of invertible elements in
$R$ with the set $R^{-1}_q$ of {\em quasi-invertible} elements,
where $u\in R^{-1}_q$ if $(1-vu)R(1-uw)=0=(1-uw)R(1-vu)$ for some
$v, w$ in $R$. (We then write $1-vu \perp 1-uw$.) It follows
easily that $u$ is a (von Neumann) regular element in $R$, and
that one may take $v=w$ and demand that $u$ and $v$ are partial
inverses for one another (i.e. $u=uvu$ and $v=vuv$). We refer to
this situation by saying that $v$ is a \emph{quasi-inverse} for
$u$. The ring $R$ is then a {\em quasi-Bass ring} (a $QB-$ring for
short) if whenever $ax+b=1$ in $R$ we can find $y$ in $R$ such
that $a+by\in R^{-1}_q$, cf. the definition of Bass rings in the
introduction to Section 2. As with the notion of stable rank, the
concept of $QB-$rings is left-right symmetric, see~\cite[Theorem
3.6]{app1}.

The theory of $QB-$rings has been developed in the
papers~\cite{app1},~\cite{app2} and~\cite{app3} with the aim of
extending as much as possible of the theory of Bass rings. The
inspiration is the corresponding series of
papers~\cite{bp2},~\cite{bp3},~\cite{bp4},~\cite{bp5}
and~\cite{bp6}, in which the theory of $C^*-$algebras of
topological stable rank one is being extended to the class of {\em
extremally rich} $C^*-$algebras, which are the $C^*-$analogues of
$QB-$rings.

We are going to prove that any surjective inverse limit of
$QB-$rings is again a $QB-$ring. To establish this in full detail
turns out to be quite intricate, and we would have preferred an
easier, more  accessible proof. The basic techniques come
from~\cite[Lemmas 1.4 and 1.5]{app3}, and we include part of the
discussion carried out there. The following Proposition examines
the situation in an easier setting, hence gives an idea of what is
going on.

In Theorem~\ref{qbmain} below we are going to use the left-handed
version of the concept of a $QB-$ring. This is done in order to
make the techniques developed in~\cite{app3} readily accessible.
Although it will not be strictly necessary for Proposition
~\ref{easycase}  we have chosen the same version there to avoid
confusion.

For the convenience of the reader we are only going to prove the
unital versions of these results. The non-unital versions follow
by straightforward, but sometimes quite exasperating computations,
where invertible and quasi-invertible elements are replaced by
{\em adversible} and {\em quasi-adversible} elements \`a la
Kaplansky, cf.~\cite[Section 4]{app1}. Basically an element $u$ in
a non-unital ring $R$ is (left/right/quasi) adversible if $1-u$ is
(left/right/quasi) invertible in some (hence any) unital ring
$\widetilde R$ containing $R$ as an ideal.

Recall that a non-unital ring $I$ has \emph{stable rank one} if
whenever $(1-x)(1-a)+b=1$ in $I^+$, where $x$, $a\in I$, we can
find $y$ in $I$ such that $1-a+yb\in (I^+)^{-1}$.

\bigskip

\begin{prop}\label{easycase}
Let $(R_n)$ be a sequence of unital $QB-$rings, and assume that we
have surjective morphisms $\pi_n:R_n\rightarrow R_{n-1}$ such that
$\ker\pi_n$ has stable rank one for every $n$. Then $R=\varprojlim
R_n$ is also a $QB-$ring.
\end{prop}

\begin{proof} Arguing along the lines of Theorem~\ref{exchange}
and Theorem~\ref{bsr}, it suffices to consider the case of a
surjective morphism $\pi\colon R\rightarrow S$, where $R$ and $S$
are unital $QB-$rings and $I=\ker\pi$ has stable rank one, along
with equations $xa+b=1$ in $R$ and $\ol{x}\,\ol{a}+\ol{b}=1$ in
$S$, so that $\pi(a)=\ol{a}$, $\pi(b)=\ol{b}$ and $\pi(x)=\ol{x}$.
Since $S$ is a $QB-$ring, there are elements $\ol{y}$ in $S$ and
$\ol{u}$ in $S_q^{-1}$ such that $\ol{a}+\ol{y}\ol{b}=\ol{u}$. We
wish to find elements $y$ in $R$ and $u$ in $R_q^{-1}$ such that
$\pi(y)=\ol{y}$, $\pi(u)=\ol{u}$ and $a+yb=u$.

First, by~\cite[Proposition 7.1]{app1} lift $\ol{u}$ to a
quasi-invertible element $u_1$ in $R$, and let $y_1$ be any lift
of $\ol{y}$. We then get
\[
a+y_1b=u_1+t\,,
\]
where $t\in\ker\pi=I$. Next, choose a quasi-inverse $v_1$ for
$u_1$. By computation
\[
x(u_1+t)+(1-xy_1)b=x(u_1+t-y_1b)+b=xa+b=1\,,
\]
so, multiplying left and right with $u_1$ and $v_1$ we have
\[
u_1x(u_1v_1+tv_1)+u_1(1-xy_1)bv_1=u_1v_1\,.
\]
Rearranging terms this gives the equation
\[
u_1x(1+tv_1)+u_1(1-xy_1)bv_1+(1-u_1x)(1-u_1v_1)=1\,.
\]
By~\cite[Lemma 4.6]{app1} we can choose $t_1$ and $s_1$ both in
$I$ such that
\[
(1+t_1)(1+tv_1)+s_1\left(u_1(1-xy_1)bv_1
+(1-u_1x)(1-u_1v_1)\right)=1\,.
\]
Since $I$ has stable rank one we can find $s_2$ in $I$ such that
\[
w=(1+tv_1)+s_2s_1\left(u_1(1-xy_1)bv_1 +(1-u_1x)(1-u_1v_1)\right)
\in (I^+)^{-1}\,,
\]
and we note that $1-w\in I$. It follows that we can write
\[
wu_1=u_1+tv_1u_1+sbv_1u_1+0=u_1+(t+sb)v_1u_1\,,
\]
where $s=s_2s_1u_1(1-xy_1)\in I$. By~\cite[Theorem 2.3]{app1} we
have that the element
\[
u_2=u_1+w^{-1}(t+sb)(1-v_1u_1)\in R^{-1}_q\,.
\]
Consequently also
\[
u=wu_2=wu_1+(t+sb)(1-v_1u_1)\in R^{-1}_q\,.
\]
Moreover, $\pi(u)=\pi(u_2) =\pi(u_1)=\ol{u}$. Define $y=y_1+s$.
Then = $\pi(y)=\ol{y}$ and
\begin{align*}
a+yb\; &=a+y_1b+sb=u_1+t+sb \\
&=u_1+(t+sb)v_1u_1+(t+sb)(1-v_1u_1)\\
&=wu_1+(t+sb)(1-v_1u_1)=u\,,
\end{align*}
as desired.
\end{proof}

\bigskip

Let $R$ be a unital ring with an ideal $I$. Consider an element
$u$ in $R_q^{-1}$ with quasi-inverse $v$, so that $u=uvu$ and
$v=vuv$, and if $p=uv$ and $q=vu$ then $(1-p)\perp (1-q)$. Let $w$
be an element in $(qRq)_q^{-1}$ such that $q-w\in I$, and choose a
quasi-inverse $w'$, satisfying the analogous equations $w=ww'w$,
$w'=w'ww'$ and $(q-ww')\perp (q-w'w)$ (whence also $q-w'\in I$).
Observe that $p-uww'v=u(q-ww')v$ is an idempotent equivalent to
$q-ww'$, whence $(p-uww'v)\perp (q-w'w)$. Set
\[
p_1=1-p, \quad p_2=p-uww'v,\quad q_1=q-w'w,\quad \text{and}\quad
q_2=1-q,
\]
and note that $p_2$, $q_1\in I$.

\bigskip

\begin{lem}\label{qblem}
Let $I$ be an ideal in a unital $QB-$ring $R$. If $u\in R^{-1}_q$
and $w\in (qRq)^{-1}_q$ as in the setup above, we consider an
element $a=uw+t_1+t_2$, where $t_i\in p_iRq_i$ for $i=1, 2$. If
$xa+b=1$ for some elements $x$ and $b$ in $R$, then there is an
element $y$ in $I$ such that $a+yb\in R_q^{-1}$.
\end{lem}

\begin{proof}
This follows from a verbatim repetition of the arguments
in~\cite[Lemma 1.4 (ii)]{app3}, using the additional assumption
that $p_iRq_i\subset I$ for $i=1$, $2$. Observe also that
by~\cite[Theorem 1.6]{app3} (cf. also~\cite[Lemma 6.1]{app1}) we
can use the following fact: either $p_i\perp q_i$, or else
$p_iRq_i$ is a $QB-$corner (for each $i$) in the sense
of~\cite[Definition 5.2]{app1}.
\end{proof}

\bigskip

\begin{theor}
\label{qbmain} If $R$ is the surjective inverse limit of a
sequence $(R_n)$ of unital $QB-$rings, then $R$ is also a
$QB-$ring.
\end{theor}

\begin{proof} As in Proposition~\ref{easycase} we only need to
consider the case of a surjective morphism $\pi:R\rightarrow S$,
where $R$ and $S$ are unital $QB-$rings, along with equations
$xa+b=1$ in $R$ and $\ol{x}\,\ol{a}+\ol{b}=1$ in $S$, so that
$\pi(a)=\ol{a}$, $\pi(b)=\ol{b}$, $\pi(x)=\ol{x}$. Since $S$ is a
$QB-$ring, there are elements $\ol{y}$ in $S$ and $\ol{u}$ in
$S^{-1}_q$ such that $\ol{a}+\ol{y}\ol{b}=\ol{u}$. We wish to find
elements $y$ in $R$ and $u$ in $R^{-1}_q$ such that
$\pi(y)=\ol{y}$, $\pi(u)=\ol{u}$ and $a+yb=u$.

Lift $\ol{u}$ to a quasi-invertible element $u_1$ in $R$ (again
by~\cite[Proposition 7.1]{app1}), and let $y_1$ be any lift of
$\ol{y}$. We then get
\[
a+y_1b=u_1+t\,,
\]
where $t\in\ker\pi=I$. By computation,
\[
x(u_1+t)+(1-xy_1)b=x(u_1+t-y_1b)+b=xa+b=1\,.
\]
If $b_1=(1-xy_1)b$ the above equation reads
\begin{equation} \label{ai}
x(u_1+t)+b_1=1
\end{equation}
Next, choose a quasi-inverse $v_1$ of $u_1$, so that
$(1-u_1v_1)\perp (1-v_1u_1)$ and $u_1=u_1v_1u_1$, $v_1=v_1u_1v_1$.
Let $p=u_1v_1$ and $q=v_1u_1$. The computations carried out
in~\cite[Lemma 1.5]{app3} can be performed in this situation also
(to the equation (\ref{ai})). Hence we obtain an element $w$ in
$(qRq)^{-1}_q\cap(q+qIq)$ and an element $z$ in $I$ such that, if
$w'$ is any quasi-inverse for $w$, and if we let $p_1=1-p$,
$p_2=p-u_1ww'v_1$, $q_1=q-w'w$ and $q_2=1-q$, then
\[
u_1+t+u_1zqb_1=w_1a_1w_2\,,
\]
where $w_1$, $w_2$ are invertible elements in $R$, and
$a_1=u_1w+t_1+t_2$, with $t_i$ in $p_iRq_i$ for $i=1, 2$.

Since $x(u_1+t)+b_1=1$, we also have
\[
x(u_1+t+u_1zqb_1)+(1-xu_1zq)b_1=1\,.
\]
Conjugating with $w_2$ we get
\[
(w_2xw_1)a_1+w_2(1-xu_1zq)b_1w_2^{-1}=1\,.
\]
By Lemma~\ref{qblem}, there is an element $y_2\in I$ such that
\[
a_1+y_2w_2(1-xu_1zq)b_1w_2^{-1}\in R_q^{-1}\,,
\]
whence $w_1a_1w_2+w_1y_2w_2(1-xu_1zq)b_1\in R_q^{-1}$. This means
that the element
\[
u=u_1+t+u_1zqb_1+w_1y_2w_2(1-xu_1zq)b_1\in R_q^{-1}\,,
\]
and also $\pi(u)=\pi(u_1)=\ol{u}$. Now let
$y=y_1+u_1zq(1-xy_1)+w_1y_2w_2(1-xu_1zq)(1-xy_1)$. Then
$\pi(y)=\ol{y}$, and
\begin{align*}
&a+yb=a+y_1b+u_1zqb_1+w_1y_2w_2(1-xu_1zq)b_1\\
=\;&u_1+t+u_1zqb_1+w_1y_2w_2(1-xu_1zq)b_1=u\,,
\end{align*}
as desired.
\end{proof}

\vskip1cm

\section{Approximate Units}

\bigskip

\noindent{\bf \boldmath$\sigma-$Unital Rings.} Contrary to popular
belief rings do not come automatically equipped with a unit. Of
course, a unit can be adjoined, passing from $R$ to $R^+$ as
decribed before, but his may destroy other desirable relations for
the ring, like being an ideal in a larger ring or being a Bass
stable ring or a $QB-$ring. Also, indiscriminate adjoining of
units may deprive the category of rings some of the applications
of general category theory.

For some rings adjoining of a unit may be unneccesary because the
ring itself already possesses elements that locally act as units.
As in ~\cite{arap} a net $(e_\lambda)_{\lambda\in\Lambda}$ in $R$
is called an {\em approximate unit} for $R$ if eventually
$e_\lambda x =xe_\lambda =x$ for every $x$ in $R$. We do not
expect that the net consists of idempotents, but in many cases it
will (for example, in the case where $R$ is an exchange ring). Far
more important is the case where we can choose $\Lambda = \mathbb
N$. We say in this case that $R$ is {\em $\sigma-$unital}.
Evidently we may then choose the approximate unit $(e_n)$ such
that $e_{n+1}e_n=e_ne_{n+1}=e_n$ for all $n$, see ~\cite[Lemma
1.5]{arap}. Note that for such a sequence to be an approximate
unit it suffices to show that for each $x$ in $R$ there is {\em
some} $n$ such that $e_nx=xe_n=x$. The same equations will then
automatically hold for all indices larger than $n$. Thus in what
follows we shall often assume that a $\sigma-$unit $(e_n)$
satisfies $e_ne_{n+1}=e_{n+1}e_n=e_n$ for all $n$.

Rings with approximate units are examples of \emph{$s-$unital}
rings (see e.g.~\cite{tomi}). If the ring is countable it is
$\sigma-$unital if and only if it is $s-$unital, see ~\cite[Lemma
2.2]{arareg}.

It should also be remarked that all $C^*-$algebras and all
classical (non-unital) Banach algebras have approximate units in
the topological sense, which implies that certain dense ideals --
representing elements with compact support or finite rank -- have
approximate units in the algebraic sense. It is instructing though
to note that the ring $\mathbb B_f(\mathcal H)$ of operators with
finite rank on an infinite dimensional Hilbert space $\mathcal H$
does not have a countable approximate unit (although it certainly
has an (uncountable) approximate unit consisting of projections).
By contrast, the ring $\mathbb M_{\infty}(R)=\varinjlim \mathbb
M_n(R)$ of finite matrices over a unital or just $\sigma-$unital
ring $R$ is $\sigma-$unital.

\bigskip

\begin{prop}\label{extsigma} Let $0\ra I\ra R \ra S\ra 0$ be a short
exact sequence of rings. If both $I$ and $S$ are $\sigma-$unital
then also $R$ is $\sigma-$unital.
\end{prop}

\begin{proof} Let $\pi\colon R \ra S$ denote the quotient map and
choose an approximate unit $(e_n)$ for $I$. Choose also a sequence
$(f_n)$ in $R$ such that $(\pi(f_n))$ is an approximate unit for
$S$.

For each $x$ in $R$ there is an $n$ such that $\pi(f_n
x)=\pi(xf_n) =\pi(x)$. Thus $(1-f_n)x\in I$ and $x(1-f_n)\in I$,
so for a suitable $m$ we have $(1-e_m)(1-f_n)x=0$ and
$x(1-f_n)(1-e_m)=0$. It follows that if we define $u_{nm}$ in $R$
by
\[
1-u_{nm}=(1-f_n)(1-e_m)(1-f_n)\,,
\]
then $u_{nm}x=xu_{nm}=x$ for some $(n,m)$ in $\mathbb N^2$. (Note
here that the use of $1$ is purely formal.)

We claim that for each $(p,q)$ in $\mathbb N^2$ there is an $n>p$
and $m>q$ such that
\[
u_{nm}u_{kl}=u_{kl}u_{nm}=u_{kl}
\]
for all $k\le p$ and $l\le q$. This is equivalent to the claim
that
\[
\big((1-f_n)(1-e_m)(1-f_n)\big)\big((1-f_k)(1-e_l)(1-f_k)\big)
=(1-f_n)(1-e_m)(1-f_n),
\]
together with the similar adjoint version. Note now that
$(1-f_k)(1-e_l)(1-f_k)=1-a_{kl}$ for some elements $a_{kl}$ in
$R$. Choosing $n$ sufficiently large we may assume that
$(1-f_n)a_{kl}\in I$ for all $k\le p$ and $l\le q$. For a
sufficiently large $m$ we therefore have $(1-e_m)(1-f_n)a_{kl}=0$
for all $k\le p$ and $l\le q$. Consequently,
\begin{align*}
&\big((1-f_n)(1-e_m)(1-f_n)\big)\big((1-f_k)(1-e_l)(1-f_k)\big)\\
=\;&(1-f_n)(1-e_m)(1-f_n)(1-a_{kl})=(1-f_n)(1-e_m)(1-f_n),
\end{align*}
as desired. The proof for the adjoint version is quite similar.

Using the claim we now inductively choose a sequence $(v_i)$ in
$R$ with $v_i=u_{n(i)m(i)}$ such that $v_iu_{nm}=u_{nm}v_i
=u_{nm}$ for all $n\le i$ and $m\le i$; but also such that
$v_{i+1}v_i=v_iv_{i+1}=v_i$ for all $i$. As we saw above, there is
for each $x$ in $R$ some $(n,m)$ in $\mathbb N^2$ such that
$u_{nm}x=xu_{nm}=x$, and it follows that if $i\ge \max\{n,m\}$
then
\[
v_ix=v_iu_{nm}x=u_{nm}x=x,
\]
and similarly $xv_i=x$. Thus $(v_i)$ is a $\sigma-$unit for $R$.
\end{proof}

\bigskip

\begin{defi}\label{def}{\rm If $I$ is an ideal in a ring $R$ we say
that $I$ {\em has a unit in $R$} if $ex=xe=x$ for some $e$ in $R$
and every $x$ in $I$. Similarly we say that $I$ {\em has a
$\sigma-$unit in $R$} if there is a sequence $e_n$ in $R$ such
that eventually $e_nx=xe_n=x$ for each $x$ in $I$.

Inspection of the proof of Proposition~\ref{extsigma} shows that
in order for $R$ to be $\sigma-$unital it suffices that $S$ is
$\sigma-$unital and $I$ has a $\sigma-$unit in $R$. In this
formulation the requirements are both necessary and sufficient.

The argument also shows that $R$ is unital if and only if $S$ is
unital and $I$ has a unit in $R$. If $e$ is a unit for $I$ in $R$
and $f$ is an element in $R$ such that $\pi(f)$ is the unit in $S$
then the element $u$ in $R$ given by $1-u=(1-f)(1-e)(1-f)$ is the
unit in $R$.}
\end{defi}

\bigskip

\begin{theor}\label{invlimsigma} let $R=\ilim R_n$ be a surjective
inverse limit of a sequence of rings with coordinate evaluations
$\rho_n\colon R\ra R_n$. Then $R$ is $\sigma-$unital if and only
if each $R_n$ is $\sigma-$unital and the ideals $\ker \rho_n$ have
a unit in $R$ eventually.
\end{theor}

\begin{proof} If $R$ is $\sigma-$unital then so is every $R_n$,
since the morphisms $\rho_n$ are surjective. The ideals $\ker
\rho_n$ form a decreasing sequence, so if they do not eventually
have a unit in $R$ none of them will. To obtain a contradiction,
let us assume this to be the case, and let $(e_n)$ be an
approximate unit for $R$.

Given $x_1$ in $R_1$ we claim that there is an $x$ in $R$ such
that $\rho_1(x)=x_1$ and either $e_2 x\ne x$ or $xe_2 \ne x$. For
if this were not the case then
\[
e_2(x+s)=(x+s)e_2=x+s
\]
for each fixed $x$ in $R$ with $\rho_1(x)=x_1$ and every $s$ in
$\ker\rho_1$. In particular, $e_2$ would be a unit for
$\ker\rho_1$ in $R$, contradicting our assumptions. Assuming
therefore, as we may, that $e_2x\ne x$ we can find a (first)
coordinate $k_2$ in which they differ. Thus, if we let
$x_2=\rho_{k_2}(x)$, then $\rho_{k_2}(e_2)x_2\ne x_2$ and
$\pi_2\circ\cdots\circ\pi_{k_2}(x_2)=x_1$. By the same argument we
can now find an element $y$ in $R$ with $\rho_{k_2}(y)=x_2$, such
that either $e_3 y\ne y$ or $ye_3 \ne y$. This produces a
coordinate number $k_3$ and an element $x_3=\rho_{k_3}(x)$ such
that either $\rho_{k_3}(e_3)x_3\ne x_3$ or $x_3\rho_{k_3}(e_3)\ne
x_3$  and $\pi_{k_2-1}\circ\cdots\pi_{k_3}(x_3)=x_2$. Continuing
by induction we find a coherent sequence $x=(x_n)$ in $(R_{k_n})$,
i.e. $\pi_{k_{n}-1} \circ \cdots \circ \pi_{k_{n+1}}(x_{n+1})=x_n$
for all $n$, such that for each $n$ either $\rho_{k_n}(e_n)x_n\ne
x_n$ or $x_n\rho_{k_n}(e_n)\ne x_n$. Identifying $x=(x_n)$ with an
element in $R$ we see that either $e_nx\ne x$ or $xe_n\ne x$ for
every $n$, contradicting the assumption that $(e_n)$ was an
approximate unit.

In the converse direction, if $R_k$ is $\sigma-$unital and
$\ker\rho_k$ has a unit in $R$, then $R$ is $\sigma-$unital by
Proposition~\ref{extsigma} and the comments after
Definition~\ref{def} applied to the extension
\[
0\ra \ker\rho_k\ra R \ra R_k\ra 0\,.
\]
\end{proof}



The Theorem above shows that in all but a few cases a surjective
inverse limit of non-unital rings will {\em not} be
$\sigma-$unital. The demand of an eventual unit for $\ker\rho_n$
in $R$ translates to the demand that with the morphisms
$\pi_m\colon R_m \ra R_{m-1}$ there is a unit $e_m$ for
$\ker\sigma_m$ in $R_m$ for each $m>n$, where
$\sigma_m=\pi_{n+1}\circ\pi_{n+2}\circ\cdots\circ\pi_{m}$, and
moreover these can be chosen coherent, i.e. $\pi_{m+1}(e_{m+1})
=e_m$ for all $m$. Clearly this is not likely to happen in actual
examples.

\vskip1cm

\section{Multiplier Rings}

\bigskip

\noindent{\bf Centralizers.} Throughout this section we assume
that $R$ is a {\em non-degenerate} ring, i.e. $xR=0$ or $Rx=0$
implies $x=0$ for every $x$ in $R$. This concept is of course only
relevant for non-unital rings, but note that rings with
approximate units will automatically be non-degenerate. Eventually
we are going to need the stronger concept that $R$ is {\em
semi-prime}, which means that $xRx=0$ implies $x=0$ for every $x$
in $R$.

Following Hochschild,~\cite[Definition 3.1]{hoch}, although using
the terminology from~\cite{johnson}, we define a {\em left
centralizer} to be a map $\lambda\colon R\rightarrow R$ such that
$\lambda(xy)= \lambda(x)y$ for all $x, y$ in R. Similarly, a {\em
right centralizer} $\rho$ is a map such that $\rho(xy)=x\rho(y)$.
A {\em double centralizer} is a pair $(\lambda,\rho)$ of
left-right centralizers satisfying the coherence relation
$x\lambda(y)=\rho(x)y$ for all $x, y$ in $R$. It is easy to see
that (in the presence of non-degeneracy) this relation alone will
imply that $\lambda$ and $\rho$ are both $R-$module maps, cf.
~\cite[Theorem 7]{johnson}. For such pairs we define a product by
setting
\[
(\lambda_1, \rho_1)(\lambda_2, \rho_2)= (\lambda_1\circ\lambda_2,
\rho_2\circ\rho_1)\,.
\]
With the obvious sum $(\lambda_1, \rho_1)+(\lambda_2, \rho_2)
=(\lambda_1 +\lambda_2, \rho_1 +\rho_2)$ the set $\mul$ of double
centralizers becomes a unital ring with the unit $1=(\id, \id)$.
We shall refer to $\mul$ as the {\em multiplier ring} for $R$.

For each $x$ in $R$ we define the left and right centralizers
$\lambda_x$ and $\rho_x$ by $\lambda_x (y)=xy$ and $\rho_x
(y)=yx$. Then $(\lambda_x, \rho_x)$ is a double centralizer, and
since $(\lambda_x,\rho_x)(\lambda_y,\rho_y)=
(\lambda_{xy},\rho_{xy})$ it follows that we have a morphism
$x\mapsto (\lambda_x,\rho_x)$ of $R$ into $\mul$. The
non-degeneracy of $R$ ensures that this morphism is injective, so
we may identify $R$ with its image in $\mul$. The easily verified
formulas
\[
\lambda\circ\lambda_x=\lambda_{\lambda(x)}, \quad
\lambda_x\circ\lambda=\lambda_{\rho(x)}, \quad
\rho\circ\rho_x=\rho_{\rho(x)}, \quad
\rho_x\circ\rho=\rho_{\lambda(x)},
\]
show that
\begin{equation}\label{mult}
(\lambda, \rho)(\lambda_x, \rho_x)=(\lambda_{\lambda(x)},
\rho_{\lambda(x)}), \quad (\lambda_x, \rho_x)(\lambda,
\rho)=(\lambda_{\rho(x)}, \rho_{\rho(x)})\,,
\end{equation}
which proves that $R$ is a two-sided ideal in $\mul$. At the same
time we see from (\ref{mult}) that $R$ is an {\em essential ideal}
in $\mul$, because $R(\lambda, \rho)=0$ or $(\lambda, \rho)R=0$
for some $(\lambda, \rho)$ in $\mul$ implies that $\lambda=0$ or
$\rho=0$, which forces $(\lambda, \rho)=0$ by the coherence
relation. It is easy to show that the multiplier ring has the
universal property that if $S$ is any ring containing $R$ as an
ideal there is a morphism $\nfi\colon S\rightarrow \mul$ extending
the embedding morphism $R\rightarrow S$. Moreover, $\nfi$ is
injective precisely when $R$ is an essential ideal in $S$. For
each $y$ in $S$ one just defines $\nfi(y) =(\lambda_y, \rho_y)$.

The multiplier ring $\mul$, or rather its quotient $\mul/R$, is an
indispensable tool in Hochschild's classification of extensions of
$R$, see the remarks at the end of the paper. It also occurs
naturally as the function-theoretic method to describe the
Stone-\v Cech compactification of a (locally compact) topological
space. In $K-$theory it is used to describe the ring $\mathbb
B(R)$ of row- and column-finite matrices over a unital ring $R$.
Specifically, if $\mathbb M_\infty(R)=\varinjlim \mathbb M_n(R)$,
identified with the ring of finite matrices over $R$, we find that
its multiplier ring is identified with $\mathbb B(R)$ (see e.g.
~\cite[Proposition 1.1]{arap}).

\begin{defi} {\rm A morphism $\pi\colon R \rightarrow S$
between rings $R$ and $S$ is called {\em proper} of $\pi(R)$ is
not contained in any proper left or right ideal of S, i.e.
$S\pi(R)=\pi(R)S=S$ (in the sense that $\pi(R)S$ denotes the set
of finite sums of products $\pi(x)y$). If $R$ is unital then $\pi$
is proper if and only if both $S$ and $\pi$ are also unital.

The origin of this notion comes from $C^*-$algebra theory, cf.
~\cite{elp} or~\cite{lor}, where one notices that if $X$ and $Y$
are locally compact Hausdorff spaces and  $R=C_0(X)$ and
$S=C_0(Y)$ denote the rings of complex-valued continuous functions
vanishing at infinity, then a morphism (which in this category
means a *-homomorphism) $\pi\colon R\rightarrow S$ is proper if
and only if it is the transposed of a {\em proper continuous map}
$f\colon Y\rightarrow X$, i.e. one for which $f^{-1}(C)$ is
compact in $Y$ for every compact subset $C$ in $X$.}

\end{defi}

\bigskip

\begin{theor}\label{multext}
Let $\pi\colon R\rightarrow S$ be a proper morphism between
non-degenerate rings $R$ and $S$. There is then a unique unital
morphism $\ol{\pi}\colon \mul\rightarrow \mathcal M(S)$ that
extends $\pi$.
\end{theor}

\begin{proof}
If $(\lambda, \rho)\in \mul$ and $y, z$ are elements in $S$ we can
write them in the form $y=\sum \pi(u_i)y_i$ and $z=\sum
z_j\pi(v_j)$ for suitable finite subsets $\{u_i\}, \{v_j\}$ in $R$
and $\{y_i\}, \{z_j\}$ in $S$. Now define
\[
\ol{\lambda}(y)=\sum \pi(\lambda(u_i))y_i\quad\text{and}\quad
\ol{\rho}(z)=\sum z_j\pi(\rho(v_j))\,.
\]
If $x\in R$ then
\[
\pi(x)\ol{\lambda}(y)=\sum \pi(x\lambda(u_i)y_i)= \sum
\pi(\rho(x)u_i)y_i=\pi(\rho(x))\sum \pi(u_i)y_i=\pi(\rho(x))y\,.
\]
Similarly, $\ol{\rho}(z)\pi(x)=z\pi(\lambda(x))$. For one thing,
this implies that the elements $\ol{\lambda}(y)$ and
$\ol{\rho}(z)$ are independent of the representations of $y$ and
$z$. If namely $\ol{\lambda'}(y)$ arises from another
representation $y=\sum \pi(u'_i)y'_i$, then
$\pi(x)(\ol{\lambda}(y)-\ol{\lambda'}(y))=\pi(\rho(x))(y-y)=0$ for
all $x$ in $R$, which forces $\ol{\lambda}(y)=\ol{\lambda'}(y)$ by
properness and non-degeneracy. On the other hand the equations
show that we have the coherence relation
\[
z\ol{\lambda}(y)=\sum z_j\pi(v_j)\ol{\lambda}(y) =\sum
z_j\pi(\rho(v_j))y=\ol{\rho}(z)y\,.
\]
We may therefore define $\ol{\pi}(\lambda, \rho)= (\ol{\lambda},
\ol{\rho})$ in $\mathcal M(S)$, and it is straightforward to check
that $\ol{\pi}$ is a unital morphism from $\mul$ to $\mathcal
M(S)$. By construction $\ol{\pi}$ is an extension of $\pi$, and
since $S$ is an essential ideal of $\mathcal M(S)$ there can be
only one extension, so $\ol{\pi}$ is unique.
\end{proof}

\bigskip

\noindent{\bf The Strict Topology.} If $R$ is a non-degenerate
ring with multiplier ring $\mul$ we can for each $a$ in $R$ define
\[
\mathcal{O}_a =\{x\in M(R) \mid xa=ax=0\}.
\]
As shown in~\cite{arap} the sets $\mathcal{O}_a$ form a subbasis
for the neighbourhood system around $0$ in the {\em strict
topology} on $\mul$. In this topology addition is continuous (by
construction) and multiplication is continuous by~\cite[Lemma
1.3]{arap}. A \emph{strict Cauchy} net $(x_\lambda)$ in $\mul$ is
a net such that each product net $(ax_\lambda)$ or $(x_\lambda a)$
is eventually constant. It follows easily from this that $\mul$ is
strictly complete, cf.~\cite[Proposition 1.6]{arap}.

If $R$ has an approximate unit $(e_\lambda)$, then
$e_\lambda\rightarrow 1$ strictly in $\mul$, and $R$ is therefore
strictly dense in $\mul$. Conversely, if $1$ is a strict limit
point for $R$ then every net $(e_\lambda)$ in $R$ converging
strictly to $1$ will be an approximate unit for $R$, and $R$ is
strictly dense in $\mul$, cf.~\cite[Proposition 1.6]{arap}.

\bigskip

\begin{prop}\label{strictcont}
If $\pi\colon R\rightarrow S$ is a proper morphism between
non-degenerate rings $R$ and $S$ then the extension
$\ol{\pi}\colon \mul\rightarrow \mathcal M(S)$ from
Theorem~\ref{multext} is strictly continuous.
\end{prop}

\begin{proof}
Let $(x_\lambda)$ be a net in $\mul$ converging strictly to $0$.
Given $y$ in $S$ we can then by properness write $y=\sum
\pi(x_i)y_i$. Consequently
\[
\ol{\pi}(x_\lambda)y = \sum \pi(x_\lambda x_i)y_i =0
\]
eventually. Similarly $y\ol{\pi}(x_\lambda)=0$ eventually, which
proves that $\ol{\pi}(x_\lambda) \rightarrow 0$ strictly in
$\mathcal M(S)$.
\end{proof}

\bigskip

The next theorem is the algebraic counterpart of a rather useful
result from $C^*-$algebra theory, known as the {\em
non-commutative Tietze extension theorem}, cf.~\cite[Proposition
3.12.10]{gkp} or~\cite[Theorem 10]{saw}. Applied to the ring
$R=C_0(X)$ of complex-valued continuous functions on some locally
compact Hausdorff space $X$, we find that $\mul =C_b(X)$, the ring
of bounded continuous functions on $X$. The content of the theorem
is then that for any bounded continuous function $f$ on a closed
subset $Y$ of $X$ (and these are the only ways quotients of the
ring $C_0(X)$ can arise), there is an extension $\ol{f}$ of $f$ to
a bounded continuous function on all of $X$. The necessity of a
{\em countable} approximate unit for $R$ reflects the fact that a
general locally compact space $X$ need not be normal (and only for
normal spaces will the Tietze extension theorem hold), but if in
addition $X$ is $\sigma-$compact, {\em then} it is normal.

\bigskip

\begin{theor}\label{tietze}
Let $\pi\colon R\rightarrow S$ be a surjective morphism between
$\sigma-$unital rings. Then the extension $\ol{\pi}\colon \mul
\rightarrow \mathcal M(S)$ from Theorem~\ref{multext} is also
surjective.
\end{theor}

\begin{proof}
Let $(e_n)$ be a countable approximate unit for $R$, and without
loss of generality assume that $e_{n+1}e_n=e_n=e_ne_{n+1}$ for all
$n$. Consider now an element $\ol{x}$ in $\mathcal M(S)$. Since
the elements $\ol{e}_n=\pi(e_n)$ form an approximate unit for $S$
we have for each $n$ that $\ol{e}_m\ol{x}\,\ol{e}_n=
\ol{x}\,\ol{e}_n$ for $m$ large enough. Passing to a subsequence
we may therefore assume that
\begin{equation}\label{norm}
\ol{e}_{n+1}\ol{x}\,\ol{e}_n = \ol{x}\,\ol{e}_n
\end{equation}
for all $n$.

Assume that for some $n\ge 2$ we have found elements $x_1, \dots,
x_n$ in $R$ with $\pi(x_k)=\ol{e}_k\ol{x}$ such that with
$e_{-1}=e_0=0$ we have
\begin{equation}\label{induct}
(x_k-x_{k-1})e_{k-3}=0=e_{k-3}(x_k-x_{k-1})
\end{equation}
for $2\le k \le n$. Choose $z$ in $R$ with $\pi(z)=
\ol{e}_{n+1}\ol{x}$. Then by (\ref{norm})
\begin{align*}
 &\pi((z-x_n)e_{n-1})=(\ol{e}_{n+1}-\ol{e}_n)\ol{x}\,\ol{e}_{n-1}=0\,,\\
 &\pi(e_{n-1}(z-x_n))=\ol{e}_{n-1}(\ol{e}_{n+1} -\ol{e}_n)\ol{x}=0\,.
\end{align*}
Thus $(z-x_n)e_{n-1}\in \ker\pi$ and  $e_{n-1}(z-x_n)\in \ker\pi$.
Now define
\[
x_{n+1}=z-(z-x_n)e_{n-1} - e_{n-1}(z-x_n)(1-e_{n-1})\,,
\]
so that $\pi(x_{n+1})=\ol{e}_{n+1}\ol{x}$. Moreover,
\begin{align*}
x_{n+1}-x_n \;&= (z-x_n)(1-e_{n-1})-e_{n-1}(z-x_n)(1-e_{n-1})\\
&=(1-e_{n-1})(z-x_n)(1-e_{n-1})\,.
\end{align*}
This means that
\[
e_{n-2}(x_{n+1}-x_n)=0=(x_{n+1}-x_n)e_{n-2}\,.
\]
By induction we can therefore find a sequence $(x_n)$ in $R$
satisfying (\ref{induct}) such that also $\pi(x_n)=\ol{e}_n\ol{x}$
for all $n\ge 3$.

For each $a$ in $R$ there is a number $k$ such that $e_ka=a=ae_k$.
It follows from (\ref{induct}) that
\[
(x_m-x_n)a=(x_m-x_n)e_ka=0\,, \quad\text{and similarly}\quad
a(x_m-x_n)=0
\]
for all $n,m \ge k+2$. Consequently $(x_n)$ is a strict Cauchy
sequence in $\mul$, hence convergent to an element $x$ in $\mul$.
Since $\ol{\pi}$ is strictly continuous by
Proposition~\ref{strictcont} we finally see that
\[
\ol{\pi}(x)=\lim \pi(x_n)=\lim \ol{e}_n\ol{x} = \ol{x}\,,
\]
as desired.
\end{proof}

\vskip1cm

\section{Multiplier Rings as Inverse Limits}

\bigskip

As before we let $R=\varprojlim R_n$ for a sequence of rings $R_n$
with surjective morphisms $\pi_n\colon R_n \rightarrow R_{n-1}$
and consider the surjective coordinate evaluation morphisms
$\rho_n\colon R\rightarrow R_n$. If $I$ is an ideal in $R$ we
evidently obtain a sequence $(I_n)$ of ideals, where
$I_n=\rho_n(R) \subset R_n$, and by restriction this defines
surjective morphisms $\pi_n\colon I_n \rightarrow I_{n-1}$.
Conversely, if we have a sequence of ideals $I_n \subset R_n$ such
that $\pi_n(I_n) =I_{n-1}$ for all $n$, we may identify the
inverse limit  $\varprojlim I_n$ with an ideal $\ol{I}$ in $R$.
Note that if we start with an ideal $I$ in $R$ then $I\subset
\ol{I}$, and in general the inclusion is strict.

\bigskip

\begin{theor} \label{multinv} If $\varprojlim R_n$ is the
surjective inverse limit of a sequence of $\sigma-$unital rings
and $I$ is an ideal of $\varprojlim R_n$ such that $\rho_n (I) =
R_n$ for every $n$ then $\mathcal M(I)=\varprojlim \mathcal
M(R_n)$.
\end{theor}

\begin{proof} Put $R=\varprojlim R_n$ and $\mathcal M
=\varprojlim \mathcal M(R_n)$. We then claim that there is a
commutative diagram
\[
\begin{CD}
I @>{\iota_0}>> R @>\rho_n>> R_{n} @>{\pi_{n}}>>
    R_{n-1}\\
@VV{\iota}V  @VV{\rho}V   @VV{\iota_{n}}V      @VV{\iota_{n-1}}V \\
\mathcal M(I) @. \mathcal M @>{\ol{\rho}_n}>> \mathcal M(R_{n})
@>{\overline\pi_n}>> \mathcal M(R_{n-1})
\end{CD}
\]

\bigskip

Here $\iota$ and $\iota_k$ for $k\ge 0$ are the natural
embeddings, and $\ol{\pi}_n$ is the surjective morphism obtained
from Theorem~\ref{tietze}. It follows that $\mathcal M$ is the
surjective inverse limit of the multiplier rings $\mathcal
M(R_n)$, and the coordinate evaluations $\ol{\rho}_n$ are
therefore also surjective.

Since the right-hand square of the diagram is commutative we can
define the morphism $\rho$ by $(x_n) \mapsto (\iota_n(x_n))$ for
every string $x=(x_n)$ in $R$, and we note that
$\ol{\rho}_n\circ\rho=\iota_n\circ\rho_n$ by this definition.
Evidently $\rho$ is injective. If moreover $y=(y_n)$ is a string
in $\mathcal M$ then for any $x$ in $R$ we have that
$\ol{\rho}_n(\rho(x)y)= \iota_n(x_n)y_n=x_ny_n$. Identifying
$(x_ny_n)$ with a string in $R$ (using the commutativity of the
diagram) we see that $\rho(R)$ is an ideal in $\mathcal M$, which
must even be essential, since $\rho(R)y=0$ implies
$R_n\ol\rho_n(y)=\rho_n(R)\ol\rho_n(y)=\ol\rho_n(\rho(R)y)=0$,
hence $\ol\rho_n(y)=0$ for all $n$ by non-degeneracy, and so
$y=0$.

We claim that $I$ is essential in $R$. For if $xI=0$ for some $x$
in $R$, then $\rho_n(x)R_n=0$ for every $n$ by our assumption on
$I$, whence $\rho_n(x)=0$ since $R_n$ is non-degenerate (being
$\sigma-$unital), and therefore $x=0$. Since the property of being
essential is hereditary it follows that $\rho(\iota_0(I))$ is an
essential ideal in $\mathcal M$. By the universal property of
multiplier rings there is therefore an injective morphism
$\nfi\colon \mathcal M\ra \mathcal M(I)$ such that
$\nfi\circ\rho\circ\iota_0=\iota$.

Each surjective morphism $\rho_n\circ\iota_0$ extends uniquely to
a (not necessarily surjective) morphism $\psi_n\colon \mathcal
M(I) \ra \mathcal M(R_n)$ by Theorem~\ref{multext}. Since
$\pi_n\circ\rho_n\circ\iota_0 =\rho_{n-1}\circ\iota_0$ and $I$ is
essential in $\mathcal M(I)$ it follows that also
$\ol\pi_n\circ\psi_n=\psi_{n-1}$ for all $n$. By the universal
property of inverse limits this means that we have a unique
morphism $\psi\colon \mathcal M(I)\ra \mathcal M$ such that
$\ol{\rho}_n\circ\psi=\psi_n$ for all $n$. It follows that
\[
\ol{\rho}_n\circ\rho\circ\iota_0=\iota_n\circ\rho_n\circ\iota_0=
\psi_n\circ\iota=\ol{\rho}_n\circ\psi\circ\iota
\]
for all $n$, which implies that $\psi\circ\iota=\rho\circ\iota_0$.

Combining these results we find that
\[
(\nfi\circ\psi)\circ\iota=\nfi\circ\rho\circ\iota_0=\iota
\quad\text{and} \quad
(\psi\circ\nfi)\circ\rho\circ\iota_0=\psi\circ\iota=\rho\circ\iota_0\,.
\]
Since $\iota(I)$ is an essential ideal in $\mathcal M(I)$ and
$\rho(\iota_0(I))$ is an essential ideal in $\mathcal M$ these
equations imply that $\nfi$ and $\psi$ are the inverse of one
another, and we have our natural isomorphism.
\end{proof}

\bigskip

\begin{corol}\label{univ} If $(R_n)$ is a sequence of
$\sigma-$unital rings with surjective morphisms $\pi_n\colon R_{n}
\rightarrow R_{n-1}$, and if $\overline\pi_n\colon \mathcal
M(R_{n})\ra \mathcal M(R_{n-1})$ denote the unique surjective
extensions of the $\pi_n$'s then
\[
\mathcal M(\varprojlim R_n) = \varprojlim \mathcal M(R_n)\,.
\]
\end{corol}

\bigskip

O'Meara proved in~\cite{om} that if $R$ is a $\sigma-$unital von
Neumann regular ring then $\mul$ is an exchange ring (but is
generally not regular). His arguments used critically the
$\sigma-$unital assumption on the ring $R$. In the next result we
show that in certain cases, multiplier rings of not necessarily
$\sigma-$unital regular rings are also exchange rings. These are
obtained as inverse limits of $\sigma-$unital regular rings, which
will generally not be $\sigma-$unital themselves, cf.
~Theorem~\ref{invlimsigma}.

\bigskip

\begin{corol}
If $R=\varprojlim R_n$ is a surjective limit of a sequence of
$\sigma-$unital von Neumann regular rings then $\mul$ is an
exchange ring.
\end{corol}

\begin{proof}
If $R_n$ is von Neumann regular and $\sigma-$unital for all $n$,
then $\mathcal M(R_n)$ is an exchange ring, by~\cite[Theorem
2]{om}. The conclusion now follows from Corollary~\ref{univ} and
Theorem~\ref{exchange}.
\end{proof}

\bigskip

Let as before $R=\varprojlim R_n$ be a surjective inverse limit of
a sequence of $\sigma-$unital rings. Assume moreover that all
$R_n$ are semi-prime. This implies that also $R$ is semi-prime.
For if $xRx=0$ for some $x$ in $R$ then $\rho_n(x)R_n\rho_n(x)=0$
for every $n$, whence $\rho_n(x)=0$ so that $x=0$. If now $I$ is
an ideal in $R$ we define the (two-sided) {\em annihilator ideal}
by
\[
I^\perp = \{ x\in R\mid xI=0\}=\{ x\in R\mid Ix=0\}\,.
\]
Then $I^\perp \cap I=0$; in fact $I^\perp$ is the largest ideal
orthogonal to $I$. In particular, $I+I^\perp$ is essential in $R$,
cf.~\cite[Exer. 6.8]{lam} or~\cite[1.1.1]{aram}.

\bigskip

\begin{defi}  {\rm Let $\varprojlim R_n$ be a surjective inverse limit
of a sequence of semi-prime rings. We say that an ideal $I$ in
$\varprojlim R_n$ is $m-${\em constant} if $I\cap \ker \rho_m =
0$. Equivalently, $I\subset (\ker\rho_m)^\perp$. We may
graphically view $\ker\rho_m$ as the strings $x=(x_n)$ in $R$ such
that $x_n=0$ for $n\le m$. Evidently these coordinate ideals
decrease, so that if $I$ is $m-$\emph{constant}, then
$I\cap\ker\rho_n=0$ for all $n\ge m$.

Since $\pi_n\circ\rho_n=\rho_{n-1}$ we see that $\ker\pi_n \subset
\rho_n(\ker\rho_m)$ for $n>m$. If therefore $I_n = \rho_n(I)$
denotes the associated sequence of ideals in $R_n$, then
$I_n\cap\ker\pi_n=0$ for $n>m$. Thus $I$ is isomorphic to $I_m$
and $I_n$ is isomorphic to $I_m$ for all $n\ge m$. In particular,
$\varprojlim I_n=I$.

Conversely, if $(I_n)$ is a sequence of ideals in $(R_n)$ such
that $\pi_n(I_{n})=I_{n-1}$ for all $n$ and $I_n\cap\ker\pi_n=0$
for all $n\ge m$ for some $m$ then $I=\varprojlim I_n$ will be an
$m-$constant ideal in $\varprojlim R_n$.

If $I$ is an $n-$constant and $J$ is an $m-$constant ideal with
$n\le m$, then $I+J\subset(\ker\rho_m)^\perp$ since
$\ker\rho_m\subset\ker\rho_n$, so $I+J$ is an $m-$constant ideal.
Thus $(\ker\rho_m)^\perp$ is the largest $m-$constant ideal. It
follows that $I_c=\bigcup (\ker\rho_m)^\perp$ is equal to the sum
of all constant ideals, and we shall refer to it as the {\em
quasi-constant ideal} of $R$.}
\end{defi}

\bigskip

The motivating example for considering constant and quasi-constant
ideals arises from the Stone-\v Cech compactification. If $X$ is a
locally compact Hausdorff space then its Stone-\v Cech
compactification $\beta X$ has the property that $C(\beta
X)=C_b(X)=\mathcal M(C_0(X))$, and this ring is always a
$C^*-$inverse limit. In the case where $X$ is also
$\sigma-$compact we can write $X=\bigcup X_n$, where each $X_n$ is
a compact subset of $X$ contained in the interior $X_{n+1}^\circ$
of $X_{n+1}$. Put $R_n =C(X_n)$ and let $\pi_n(f)= f| X_{n-1 }$
for each $f$ in $C(X_n)$. Then $C_b(X) = \varprojlim R_n$. The
large constant ideals will be of the form $(\ker\pi_m)^\perp =
C_0(X_m^\circ)$, so the quasi-constant ideal of $C_b(X)$ can (in
this category) be identified with the norm-closed ideal $C_0(X)$.

If the construction is carried out in a purely algebraic context
we find that the ring $C(X)$ of all continuous functions on $X$ is
the (algebraic) inverse limit of the sequence $(C(X_n))$, and the
quasi-constant ideal of $C(X)$ is the ring $C_c(X)$ of continuous
functions with compact supports.

Theorem~\ref{multinv} provides an immediate generalization of this
construction:

\begin{corol} If $R = \varprojlim R_n$ is the surjective
inverse limit of a sequence of $\sigma-$unital rings such that the
quasi-constant ideal $ I_c$ of $R$ satisfies $\rho_n (I_c) = R_n$
for every $n$, then $\mathcal M(I_c) =\varprojlim \mathcal
M(R_n)$.
\end{corol}

\bigskip

\begin{theor} \label{idealinv} Let $(I_n)$ and $(J_n)$ be
two sequences of ideals in a semi-prime $\sigma-$unital ring $R$,
one increasing, the other decreasing, but such that $I_n\cap
J_n=0$ for all $n$. If $I = \bigcup I_n$ is essential in $R$ and
$I + J_n = R$ for every $n$, then with $R_n = R/J_n$ and
$\pi_n\colon R_{n}\to R_{n-1}$ the natural morphisms we have an
embedding $\ol\nfi\colon R\rightarrow \varprojlim R_n$ such that
$\ol\nfi (I)$ is an ideal. Moreover, $\mathcal M(I) =\varprojlim
\mathcal M(R_n)$.
\end{theor}

\begin{proof} If $x\in \bigcap J_n$ then it annihilates
$I_n$ for every $n$, whence $x\in I^\perp$. But then $x=0$ since
$I$ is essential and $R$ is semi-prime. Thus our assumptions imply
that $\bigcap J_n = 0$.

The quotient morphisms $\nfi_n\colon R\to R_n$ satisfy
$\pi_n\circ\nfi_n=\nfi_{n-1}$ for all $n$, and therefore define a
morphism $\ol\nfi\colon R\to \varprojlim R_n$ such that
$\rho_n\circ\ol\varphi=\varphi_n$ for all $n$, where
$\rho_n:\varprojlim R_n\rightarrow R_n$ are the coordinate
evaluations. Since $\ker\ol\nfi =\bigcap J_n = 0$, this is an
embedding.

Observe that $\nfi_n(I_m)=(I_m+J_n)/J_n$ is an ideal in $R_n$ for
every $n$ and $m$. Since $\ker\pi_n=J_{n-1}/J_n$ we see moreover
that $\nfi_n(I_m)\cap\ker\pi_n=0$ for $n>m$. Consequently
$(\nfi_n(I_m))$ is a coherent sequence of ideals in $(R_n)$, all
isomorphic for $n\ge m$, thus giving rise to the $m-$constant
ideal $\ol\nfi(I_m)$ in $\varprojlim R_n$.

It follows from this that $\ol\nfi(I)$ is an ideal in $\varprojlim
R_n$, isomorphic to $I$ (and contained in the quasi-constant ideal
of $\varprojlim R_n$). Since by assumption
\[
\rho_n(\ol\nfi(I))=\nfi_n(I)=(I+J_n)/J_n=R_n\,,
\]
it follows from Theorem~\ref{multinv} that $\mathcal M(I)=
\varprojlim \mathcal M(R_n)$.
\end{proof}

\bigskip

\begin{corol}\label{duideal}
Let $(I_n)$ and $(J_n)$ be two sequences of ideals in a semi-prime
$\sigma-$unital ring $R$, one increasing, the other decreasing,
but such that $I_n\cap J_n=0$ for every $n$. If $\bigcup I_n =R$
and each quotient $R_n = R/J_n$ is unital, then with $\pi_n\colon
R_n\to R_{n-1}$ the natural morphisms we have an embedding of $R$
as an ideal in $\varprojlim R_n$ such that $\mul =\varprojlim R_n
$.
\end{corol}

\bigskip

\begin{exem}
{\rm Let $(e_n)$ be a countable approximate unit consisting of
central elements in a semi-prime ring $R$. Assuming, as we may,
that $e_ne_{n+1}=e_n$ for all $n$ we consider the principal ideals
$I_n=e_nR$ and $J_n=(1-e_{n+1})R$ in $R$. Note that $e_n\in I_n$,
whereas $J_n$ consists of elements of the form $x-e_{n+1}x$, where
$x\in R$. Evidently $I_n\cap J_n=0$ and $\bigcup I_n=R$. Moreover
$R_n=R/J_n$ is unital with unit $e_{n+1}+J_n$. It follows from
Corollary~\ref{duideal} that $\mul= \varprojlim R_n$.

Along the same lines we may assume that the semi-prime ring $R$
has a countable approximate unit $(e_n)$, such that
$e_{n+1}x=x=xe_{n+1}$ for every $x$ in the principal ideal $I_n$
generated by $e_n$. This could be accomplished by placing suitable
finiteness conditions on $I_n$, for example finite dimensionality.
As before we let $J_n$ be the principal ideal in $R$ generated by
$1-e_{n+1}$ (consisting of sums of elements of the form
$x-xe_{n+1},\; y-e_{n+1}y$ and $xy-xe_{n+1}y$, with $x, y$ in
$R$). By assumption $I_n\cap J_n=0$ and $\bigcup I_n=R$. Moreover,
each quotient $R_n=R/J_n$ is unital with unit $e_{n+1}+J_n$, so
again $\mul =\varprojlim R_n$.

Note that in both examples the quotients $R_{n-1}=R/J_{n-1}$ can
be identified with the local rings $R_{e_n}$ ($= e_nRe_n$ as an
additive group, but with the product
$(e_nxe_n)(e_nye_n)=e_nxe_nye_n$).}
\end{exem}

\bigskip

\noindent{\bf Pullbacks.} Since the functor $\varprojlim$ has a
left adjoint it commutes with extensions and pullbacks. Thus if
for each $n$ we have a commutative diagram of rings with morphisms
\[
\begin{CD}
R_n  @>{\alpha_n}>>               T_n  @<{\beta_n}<<
   S_n\\
@VV{\pi_n}V                        @VV{\tau_n}V     @VV{\rho_n}V\\
R_{n-1} @>{\alpha_{n-1}}>>   T_{n-1}   @<{\beta_{n-1}}<< S_{n-1}
\end{CD}
\]
we can form the pullback rings and obtain induced morphisms
\[
\begin{CD}
R_n\oplus_{T_n} S_n @>{\sigma_n}>> R_{n-1}\oplus_{T_{n-1}}
S_{n-1}\,.
\end{CD}
\]
If for simplicity we put
\[
R=\varprojlim R_n, \quad S=\varprojlim S_n \quad \text{and}\quad
T=\varprojlim T_n
\]
we also have induced morphisms $R\rightarrow T \leftarrow S$, and
as in~\cite[Proposition 4.16]{pulp} we obtain a natural
isomorphism
\[
\varprojlim  \left(R_n \oplus_{T_n} S_n \right)  = R\oplus_T S
\]

\medskip

\noindent{\bf Extensions.} If we have a commutative diagram of
extensions
\begin{equation}\label{diagram}
\begin{CD}I_n  @>{\iota_n}>>      R_n  @>{\alpha_n}>>
    S_n\\
@VV{\pi_n}V           @VV{\omega_n}V     @VV{\theta_n}V\\
I_{n-1}  @>{\iota_{n-1}}>>   R_{n-1}   @>{\alpha_{n-1}}>> S_{n-1}
\end{CD}
\end{equation}
then, again by category theory, we obtain an extension
\begin{equation}\label{limext}
\begin{CD}
\varprojlim I_n @>{\ol\iota}>>\varprojlim R_n
@>{\ol\alpha}>>\varprojlim S_n\,.
\end{CD}
\end{equation}

We know from~\cite{hoch} that every extension is determined (up to
an obvious  isomorphism) by the morphism $\sigma$ (the {\em
Hochschild invariant}) occurring in the natural diagram below:
\[
\begin{CD}
I  @>{\iota}>>      R @>{\alpha}>>              S\\
@|           @VV{\nfi}V     @VV{\sigma}V\\
I  @>{\iota}>>   \mathcal M(I)   @>{\beta}>>  \mathcal Q(I)
\end{CD}
\]
Here $\mathcal M(I)$ is the multiplier ring of the ideal $I$,
assumed to be non-degenerate, and $\mathcal Q(I)$, the corona
ring, is the quotient $\mathcal M(I)/I$. The morphism $\nfi$ is
the one obtained by the universal property of $\mathcal M(I)$, and
$\sigma$ is the morphism induced by $\nfi$ to make the whole
diagram commutative. The salient fact to know is that the
right-hand square in the diagram is a pullback.

As a consequence of our results about multipliers of inverse
limits, in particular Corollary~\ref{univ}, we can compute the
Hochschild invariant of an inverse limit of extensions:

\bigskip

\begin{prop}
Given a coherent sequence of extensions as in the diagram
(\ref{diagram}), each determined by a Hochschild invariant
$\sigma_n\colon S_n \rightarrow \mathcal Q(I_n)$, where $(I_n)$ is
sequence of $\sigma-$unital ideals and all the morphisms $\pi_n$
and $\omega_n$ (hence also $\theta_n$) are surjective, then since
$\mathcal Q(\varprojlim I_n)=\varprojlim \mathcal Q(I_n)$, the
Hochschild invariant for the extension (\ref{limext}) is the
induced morphism $\ol\sigma$ occurring in the diagram:
\[
\begin{CD}
\varprojlim I_n @>{\ol\iota}>> \varprojlim R_n @>{\ol\nfi}>>
\varprojlim S_n\\
@|           @VV{\ol\nfi}V     @VV{\ol\sigma}V\\
\varprojlim I_n @>{\ol\iota}>>\varprojlim \mathcal M(I_n)
@>{\ol\psi}>>\varprojlim \mathcal Q(I_n)
\end{CD}
\]
\end{prop}

\end{document}